\documentclass[12pt,thmsa]{article}
\usepackage{amsfonts}
\usepackage{amsmath}
\usepackage{amssymb}

\input{tcilatex}
\begin{document}

\title{Abstract Harmonic Analysis on the General Linear Group $GL(n,\mathbb{R%
})$  }
\author{Kahar El-Hussein \\
\textit{Department of Mathematics, Faculty of Science, }\\
\textit{\ Al Furat University, Dear El Zore, Syri\textit{a}\ and}\\
\textit{Department of Mathematics, Faculty of Arts Science Al Quryyat, }\\
\textit{\ Al-Jouf University, KSA }\\
\textit{E-mail : kumath@ju.edu.sa, kumath@hotmail.com }}
\maketitle

\begin{abstract}
Let $GL(n,\mathbb{R})$ be the general linear group and let $SL(n,\mathbb{R}%
)=KAN$ be the Iwasawa decomposition of real connected semisimple Lie group $%
SL(n,\mathbb{R})$. We adopt the technique in my paper $[12]$ to generalize
the definition of the Fourier transform and to obtain the Plancherel theorem
for the group $SL(n,\mathbb{R}).$ Since $GL_{+}(n,\mathbb{R})$ $=SL(n,%
\mathbb{R})\times \mathbb{R}_{+}^{\ast }$ is the direct product of $SL(n,%
\mathbb{R})$ and $\mathbb{R}_{+}^{\ast }$, so we obtain the Plancherel
formula for theorem $GL_{+}(n,\mathbb{R})$. In the end we prove that $%
GL_{-}(n,\mathbb{R})$ has a structure of group, which is isomorphic onto $%
GL_{+}(n,\mathbb{R}),$ and leads us to obtain the Plancherel theorem for $%
GL(n,\mathbb{R}).$
\end{abstract}



\section{\textbf{\ Introduction.}}

\bigskip\ \textbf{1.1. }As a manifold, $GL(n,\mathbb{R})$ is not connected
but rather has two connected components, the matrices with positive
determinant and the ones with negative determinant. The identity component,
denoted by $GL_{+}(n,\mathbb{R})$, consists of the real $n\times n$ matrices
with positive determinant. This is also a Lie group of dimension $n^{2}$; it
has the same Lie algebra as $GL(n,\mathbb{R})$. The group $GL(n,\mathbb{R})$
is also noncompact. The maximal compact subgroup of $GL(n,\mathbb{R})$ is
the orthogonal group $O(n)$, while the maximal compact subgroup of $GL_{+}(n,%
\mathbb{R})$ is the special orthogonal group $SO(n)$. As for $SO(n)$, the
group $GL_{+}(n,\mathbb{R})$ is not simply connected (except when $n=1$),
but rather has a fundamental group isomorphic to $%
\mathbb{Z}
$ for $n=2$ or $%
\mathbb{Z}
_{2}$ for $n>2.$ Various physical systems, such as crystals and the hydrogen
atom, can be modelled by symmetry groups. In quantum mechanics, we know that
a change of frame a gauge transform leaves the probability of an outcome
measurement invariant (well, the square modulus of the wave-function, i.e.
the probability), because it is just a multiplication by a phase term
theories of gravitation based, respectively, on the general linear group $%
GL(n,\mathbb{R})$ and its inhomogeneous extension $IGL(n,\mathbb{R})$. The
topological groups $GL(n,\mathbb{R})$ and $GL(n,\mathbb{%
\mathbb{C}
})$ have an inexhaustibly rich structure and importance in all parts of
modern mathematics: analysis, geometry, topology, representation theory,
number theory, ect.... . The general linear group $GL(n,\mathbb{R})$ is
decomposed into a Markov-type Lie group and an abelian scale group. The
problem of finding the explicit Plancherel formulas for semisimple Lie
groups has been solved completely in the case of complex semisimple Lie
groups $[13]$. Moreover Harish-Chandra showed $[15]$ that the problem is
solved also for a real semisimple Lie group having only one conjugate class
of Cartan subgroups. In the case of real semisimple Lie groups with several
conjugate classes of Cartan subgroups, the problem is very difficult to
attack. The problem was taken up and solved for $SL(2,\mathbb{R})$ by V.
Bargman, $[1]$, Harish-Chandra $[14]$ and L. Pukanszky $[21]$ also for the
universal covering group of $SL(2,\mathbb{R})$ by L. Pukanszky. In this
paper we will establish the Plancherel formula for the general linear group
on $GL(n,\mathbb{R}).$ The method is based on my paper $[$ $12]$

\section{\protect\bigskip Notation and Results}

\textbf{2.1. }The fine structure of the nilpotent Lie groups will help us to
do the Fourier transform on a simply connected nilpotent Lie groups $N.$ As
well known any group connected and simply connected $N$ has the following
form%
\begin{equation}
N=\left( 
\begin{array}{cccccccccccc}
1 & x_{1}^{1} & x_{1}^{2} & x_{1}^{3} & . & . & . & . & . & x_{1}^{n-2} & 
x_{1}^{n-1} & x_{1}^{n} \\ 
0 & 1 & x_{2}^{2} & x_{2}^{3} & . & . & . & . & . & x_{2}^{n-2} & x_{2}^{n-1}
& x_{2}^{n} \\ 
0 & 0 & 1 & x_{3}^{3} & . & . & . & . & . & x_{3}^{n-2} & x_{3}^{n-1} & 
x_{3}^{n} \\ 
0 & 0 & 0 & 1 & . & . & . & . & . & x_{4}^{n-2} & x_{4}^{n-1} & x_{4}^{n} \\ 
. & . & . & . & . & . & . & . & . & . & . & . \\ 
. & . & . & . & . & . & . & . & . & . & . & . \\ 
. & . & . & . & . & . & . & . & . & . & . & . \\ 
. & . & . & . & . & . & . & . & . & . & . & . \\ 
. & . & . & . & . & . & . & . & . & x_{n-2}^{n-2}. & .x_{n-2}^{n-1} & 
x_{n-2}^{n} \\ 
0 & 0 & 0 & 0 & . & . & . & . & . & 1 & x_{n-1}^{n-2} & x_{n=1}^{n} \\ 
0 & 0 & 0 & 0 & . & . & . & . & . & 0 & 1 & x_{n}^{n} \\ 
0 & 0 & 0 & 0 & . & . & . & . & . & 0 & 0 & 1%
\end{array}%
\right)
\end{equation}

\bigskip As shown, this matrix is formed by the subgroup $\mathbb{R}$, $%
\mathbb{R}^{2}$,...., $\mathbb{R}^{n-1}$, and $\mathbb{R}^{n}$ 
\begin{equation}
\left( \mathbb{R=}\left[ 
\begin{array}{c}
x_{1}^{1} \\ 
1 \\ 
0 \\ 
0 \\ 
. \\ 
. \\ 
. \\ 
. \\ 
. \\ 
0 \\ 
0 \\ 
0%
\end{array}%
\right] ,\mathbb{R}^{2}=\left[ 
\begin{array}{c}
x_{1}^{2} \\ 
x_{2}^{2} \\ 
1 \\ 
0 \\ 
. \\ 
. \\ 
. \\ 
. \\ 
. \\ 
0 \\ 
0 \\ 
0%
\end{array}%
\right] ,..,\mathbb{R}^{n-1}=\left[ 
\begin{array}{c}
x_{1}^{n-1} \\ 
x_{2}^{n-1} \\ 
x_{3}^{n-1} \\ 
x_{4}^{n-1} \\ 
. \\ 
. \\ 
. \\ 
. \\ 
.x_{n-2}^{n-1} \\ 
x_{n-1}^{n-2} \\ 
1 \\ 
0%
\end{array}%
\right] ,\mathbb{R}^{n}=\left[ 
\begin{array}{c}
x_{1}^{n} \\ 
x_{2}^{n} \\ 
x_{3}^{n} \\ 
x_{4}^{n} \\ 
. \\ 
. \\ 
. \\ 
. \\ 
x_{n-2}^{n} \\ 
x_{n=1}^{n} \\ 
x_{n}^{n} \\ 
1%
\end{array}%
\right] \right)
\end{equation}%
$\ \ \ \ \ \ \ \ \ \ \ \ \ \ \ \ \ \ \ \ $

Each $\mathbb{R}^{i}$ is a subgroup of $N$ of dimension $i$ , $1\leq i\leq
n, $ put $d=n+(n-1)+....+2+1=\frac{n(n+1)}{2},$ which is the dimension of $N$
. According to $[9],$ the group $N$ is isomorphic onto the following group

\begin{equation}
(((((\mathbb{R}^{n}\rtimes _{\rho _{n}})\mathbb{R}^{n-1})\rtimes _{\rho
_{n-1}})\mathbb{R}^{n-2}\rtimes _{\rho _{n-2}}.....)\rtimes _{\rho _{2}}%
\mathbb{R}^{2})\rtimes _{\rho _{1}}\mathbb{R}
\end{equation}

\bigskip That means

\begin{equation}
N\simeq (((((\mathbb{R}^{n}\rtimes _{\rho _{n}})\mathbb{R}^{n-1})\rtimes
_{\rho _{n-1}})\mathbb{R}^{n-2}\rtimes _{\rho _{n-2}}.....)\rtimes _{\rho
_{4}}\mathbb{R}^{3}\rtimes _{\rho _{3}}\mathbb{R}^{2})\rtimes _{\rho _{2}}%
\mathbb{R}
\end{equation}

\textbf{2.2}$.$ Denote by $L^{1}(N)$ the Banach algebra that consists of all
complex valued functions on the group $N$, which are integrable with respect
to the Haar measure of $N$ and multiplication is defined by convolution on $%
N $ as follows:%
\begin{equation}
g\ast f(X)=\int\limits_{N}f(Y^{-1}X)g(Y)dY
\end{equation}%
for any $f\in L^{1}(N)$ and $g\in L^{1}(N),$ where $X=(X^{1},$ $X^{2},$ $%
X^{3},....,X^{n-2},X^{n-1},X^{n}),$ $%
Y=(Y^{1},Y^{2},Y^{3},....,Y^{n-2},Y^{n-1},Y^{n}),$ $X^{1}=x_{1}^{1},$ $%
X^{2}=(x_{1}^{2},x_{2}^{2}),$ $X^{3}=(x_{1}^{3},x_{2}^{3},x_{3}^{3})$ $%
,...., $ $%
X^{n-2}=(x_{1}^{n-2},x_{2}^{n-2},x_{3}^{n-2},x_{4}^{n-2},...,x_{n-2}^{n-2}),$
$%
X^{n-1}=(x_{1}^{n-1},x_{2}^{n-1},x_{3}^{n-1},x_{4}^{n-1},...,x_{n-2}^{n-1},x_{n-1}^{n-1}), 
$ $%
X^{n}=(x_{1}^{n},x_{2}^{n},x_{3}^{n},x_{4}^{n},...,x_{n-2}^{n},x_{n-1}^{n},x_{n}^{n}) 
$ and $dY=dY^{1}dY^{2}dY^{3},....,dY^{n-2}dY^{n-1}dY^{n}$ is the Haar
measure on $N$ and $\ast $ denotes the convolution product on $N.$ We denote
by $L^{2}(N)$ its Hilbert space. Let $\mathcal{U}\;$be the complexified
universal enveloping algebra of the real Lie algebra $\underline{g}$\ of $N$%
; which is canonically isomorphic to the algebra of all distributions on $N$
supported by $\left\{ 0\right\} ,$ where $0$ is the identity element of $N$.
For any $u\in \mathcal{U}$ one can define a differential operator $P_{u}$ on 
$N$ as follows:%
\begin{equation}
P_{u}f(X)=u\ast f(X)=\int\limits_{N}f(Y^{-1}X)u(Y)dY
\end{equation}

\bigskip The mapping $u\rightarrow P_{u}$ is an algebra isomorphism of $%
\mathcal{U}$ onto the algebra of all invariant differential operators on $N$

\section{Fourier Transform and Plancherel Theorem on $N.$}

\bigskip \textbf{3.1. }For each $\ 2\leq i\leq n,$ let 
\begin{eqnarray*}
&&K_{i-1} \\
&=&\left\{ \underbrace{\left\{ \underbrace{\left\{ \underbrace{\left\{ 
\underbrace{\mathbb{R}^{n}}_{\rtimes _{\rho _{n}}}\right\} \times \mathbb{R}%
^{n-1}\times \mathbb{R}^{n-1}}_{\rtimes _{_{\rho _{n-1}}}}\right\} \times 
\mathbb{R}^{n-2}\times \mathbb{R}^{n-2}\times .....}_{\rtimes _{\rho
i+1}}\right\} \times \mathbb{R}^{i}\times \mathbb{R}^{i}}_{\rtimes _{\rho
i}}\right\} \times \mathbb{R}^{i-1}\mathbb{\rtimes }_{\rho _{i}}\mathbb{R}%
^{i-1}
\end{eqnarray*}%
be the group with the following law

\begin{eqnarray}
&&(X....,X^{i},X^{i},U^{i-1},X^{i-1}).(Y....,Y^{i},Y^{i},V^{i-1},Y^{i-1}) \\
&=&[(X....,X^{i},X^{i}).(\rho
_{i}(X^{i})(Y....,Y^{i},Y^{i}))],U^{i-1}+V^{i-1},X^{i-1}+Y^{i-1})  \notag
\end{eqnarray}

For any $X\in L_{i+1},Y\in L_{i+1},$ $X^{i}\in \mathbb{R}^{i},U^{i-1}\in 
\mathbb{R}^{i-1},X^{i-1}\in \mathbb{R}^{i-1},Y^{i}\in \mathbb{R}%
^{i},V^{i-1}\in \mathbb{R}^{i-1},Y^{i-1}\in \mathbb{R}^{i-1}.$ So we get for 
$i=2$

\ 
\begin{eqnarray*}
&&K_{1} \\
&=&\left\{ \underbrace{\left\{ \underbrace{\left\{ \underbrace{\left\{ 
\underbrace{\mathbb{R}^{n}}_{\rtimes _{\rho _{n}}}\right\} \times \mathbb{R}%
^{n-1,1}\rtimes _{\rho _{n-1}}\mathbb{R}^{n-1}}_{\rtimes _{_{\rho
_{n-1}}}}\right\} \times \mathbb{R}^{n-2,1}\rtimes _{\rho _{3}}\mathbb{R}%
^{n-2}\times .....}_{\rtimes _{\rho 3}}\right\} \times \mathbb{R}%
^{2,1}\times \mathbb{R}^{2}}_{\rtimes _{\rho _{2}}}\right\} \times \mathbb{R}%
^{1,1}\mathbb{\rtimes }_{\rho _{2}}\mathbb{R}^{1}
\end{eqnarray*}%
and for $i=n$ , we get 
\begin{equation}
K_{n-1}=\left\{ \underbrace{\mathbb{R}^{n}}_{\rtimes _{\rho _{n}}}\right\}
\times \mathbb{R}^{n-1}\rtimes _{\rho _{n}}\mathbb{R}^{n-1}  \notag
\end{equation}

\bigskip \textbf{3.2. }Let $M=\mathbb{R}^{n}\times \mathbb{R}^{n-1,1}\times 
\mathbb{R}^{n-2,1}\times .....\times \mathbb{R}^{3,1}\times \mathbb{R}%
^{2,1}\times \mathbb{R}^{1,1}\mathbb{=R}^{d}$ be the Lie group, which is the
direct product of $\mathbb{R}^{n,1}=\mathbb{R}^{n},\mathbb{R}^{n-1,1}=%
\mathbb{R}^{n-1},\mathbb{R}^{n-2,1}=\mathbb{R}^{n-2},.....,\mathbb{R}^{3,1}=%
\mathbb{R}^{3},\mathbb{R}^{2,1}=\mathbb{R}^{2}$ and $\mathbb{R}^{1,1}=%
\mathbb{R}^{1}$. Denote by $L^{1}(M)$ the Banach algebra that consists of
all complex valued functions on the group $M$, which are integrable with
respect to the Lebesgue measure on $M$ and multiplication is defined by
convolution on $M$ as: 
\begin{equation}
g\ast _{c}f(X)=\int\limits_{M}f(X-Y)g(Y)dY
\end{equation}%
for any $f\in L^{1}(M),$ $g\in L^{1}(M),$ where $\ast _{c}$ signifies the
convolution product on the abelian group $M.$ We denote again by $\mathcal{U}%
\;$be the complexified universal enveloping algebra of the real Lie algebra $%
\underline{m}$\ of $M$; which is canonically isomorphic to the algebra of
all distributions on $M$ supported by $\left\{ 0\right\} ,$ where $0$ is the
identity element of $M$. For any $u\in \mathcal{U}$ one can define a
differential operator $Q_{u}$ on $M$ as follows:%
\begin{equation}
Q_{u}f(X)=u\ast _{c}f(X)=f\ast _{c}u(X)=\int\limits_{M}f(X-Y)u(Y)dY
\end{equation}

\bigskip The mapping $u\rightarrow Q_{u}$ is an algebra isomorphism of $%
\mathcal{U}$ onto the algebra of all invariant differential operators (the
algebra of all linear differential operators with constant coefficients) on $%
M$

The group $N$ can be identified with the subgroup%
\begin{equation}
N=\left\{ \underbrace{\left\{ \underbrace{\left\{ \underbrace{\left\{ 
\underbrace{\mathbb{R}^{n}}_{\rtimes _{\rho _{n}}}\right\} \times \mathbb{\{}%
0\mathbb{\}}\rtimes _{\rho _{n-1}}\mathbb{R}^{n-1}}_{\rtimes _{_{\rho
_{n-1}}}}\right\} \times \mathbb{\{}0\mathbb{\}}\rtimes _{\rho _{3}}\mathbb{R%
}^{n-2}\times .....}_{\rtimes _{\rho 3}}\right\} \times \mathbb{\{}0\mathbb{%
\}}\times \mathbb{R}^{2}}_{\rtimes _{\rho _{2}}}\right\} \times \mathbb{\{}0%
\mathbb{\}\rtimes }_{\rho _{2}}\mathbb{R}
\end{equation}%
of $K_{1}$ and $M$ can be identified with the subgroup 
\begin{eqnarray*}
&&M \\
&=&\left\{ \underbrace{\left\{ \underbrace{\left\{ \underbrace{\left\{ 
\underbrace{\mathbb{R}^{n}}_{\rtimes _{\rho _{n}}}\right\} \times \mathbb{R}%
^{n-1,1}\mathbb{\rtimes }_{\rho _{n}}\mathbb{\{}0\mathbb{\}}}_{\rtimes
_{_{\rho _{n-1}}}}\right\} \times \mathbb{R}^{n-2,1}\mathbb{\rtimes }_{\rho
_{n-1}}\mathbb{\{}0\mathbb{\}}\times .....}_{\rtimes _{\rho 3}}\right\}
\times \mathbb{R}^{2,1}\mathbb{\rtimes }_{\rho _{3}}\mathbb{\{}0\mathbb{\}}}%
_{\rtimes _{\rho _{2}}}\right\} \times \mathbb{R}^{1,1}\mathbb{\rtimes }%
_{\rho _{2}}\mathbb{\{}0\mathbb{\}}
\end{eqnarray*}

In this paper, we show how the Fourier transform on the vector group $%
\mathbb{R}^{d}$ can be generalized on $N$ and obtain the Plancherel theorem.

\textbf{Definition 3.1.} \textit{For} $1\leq i\leq n,$ \textit{let} $%
\mathcal{F}^{i}$\textit{\ be the Fourier transform on} $\mathbb{R}^{i}$ 
\textit{and }$0\leq j\leq n-1,$ \textit{let }$\dprod\limits_{0\leq l\leq j}%
\mathbb{R}^{n-l}=(..((((\mathbb{R}^{n}\rtimes _{\rho _{n}})\mathbb{R}%
^{n-1})\rtimes _{\rho _{n-1}})\mathbb{R}^{n-2}\rtimes _{\rho
_{n-2}})..\times _{\rho _{n-j}}\mathbb{R}^{n-j}),$ \textit{and let }$%
\dprod\limits_{0\leq l\leq j}\mathcal{F}^{n-l}$ $=\mathcal{F}^{n}\mathcal{F}%
^{n-l}\mathcal{F}^{n-2}....\mathcal{F}^{n-j},$\textit{we can define the
Fourier transform on } $\dprod\limits_{0\leq l\leq n-1}\mathbb{R}^{n-l}=%
\mathbb{R}^{n}\rtimes _{\rho _{n}}\mathbb{R}^{n-1}$ $\rtimes _{\rho _{n-1}}%
\mathbb{R}^{n-2}\rtimes _{\rho _{n-2}}.....\rtimes _{\rho _{3}}\mathbb{R}%
^{2}\rtimes _{\rho _{2}}\mathbb{R}^{1}$\textit{as}%
\begin{eqnarray}
&&\mathcal{F}^{n}\mathcal{F}^{n-1}\mathcal{F}^{n-2}....\mathcal{F}^{2}%
\mathcal{F}^{1}f(\lambda ^{n},\text{ }\lambda ^{n-1},\lambda ^{n-2},\text{%
.......,}\lambda ^{2},\lambda ^{1})  \notag \\
&=&\dint\limits_{N}f(X^{n},X^{n-1},...,X^{2},X^{1})e^{-\text{ }i\langle 
\text{ }(\lambda ^{n},\text{ }\lambda ^{n-1}),(X^{n},X^{n-1}),..,(\lambda
^{2},\text{ }\lambda ^{1}),(X^{2},X^{1})\rangle }\text{ }  \notag \\
&&dX^{n}dX^{n-1}...dX^{2}dX^{1}
\end{eqnarray}%
\textit{for any} $f\in L^{1}(N),$ \textit{where }$%
X=(X^{n},X^{n-1},..,X^{n},X^{n-1}),$ $\mathcal{F}^{d}=\mathcal{F}^{n}%
\mathcal{F}^{n-1}\mathcal{F}^{n-2}....\mathcal{F}^{2}\mathcal{F}^{1}$\ 
\textit{is} \textit{the classical Fourier transform on} $N,$ $%
dX=dX^{n}dX^{n-1}...dX^{2}dX^{1},$ $\lambda =(\lambda ^{n},$ $\lambda
^{n-1},\lambda ^{n-2},..$,$\lambda ^{2},\lambda ^{1}),$ \textit{and}

$\langle (\lambda ^{n},$ $\lambda ^{n-1}),(X^{n},X^{n-1}),..,(\lambda ^{2},$ 
$\lambda ^{1}),(X^{2},X^{1})\rangle =\dsum\limits_{i=1}^{n}X_{i}^{n}\lambda
_{i}^{n}+\dsum\limits_{j=1}^{n-1}X_{j}^{n-1}\lambda
_{j}^{n-1}+..+\dsum\limits_{i=1}^{2}X_{i}^{2}\lambda _{i}^{2}+X^{1}\lambda
^{1}$

\bigskip \textbf{Plancherels theorem 3.2}. For any function $f\in L^{1}(N),$
we have%
\begin{eqnarray}
&&\dint_{N}\left\vert \mathcal{F}^{d}f(\xi
^{n},X^{n-1},X^{n-2},..,X^{2},X^{1})\right\vert
^{2}dXdX^{n-1}dX^{n-2}..dX^{2}dX^{1}  \notag \\
&&\dint_{N}\left\vert \mathcal{F}^{d}f(\xi ^{n},\lambda ^{n-1},\lambda
^{n-2},..,\lambda ^{2},\lambda ^{1})\right\vert ^{2}d\xi ^{n}d\lambda
^{n-1}d\lambda ^{n-2}..d\lambda ^{2}d\lambda ^{1}
\end{eqnarray}

\bigskip \textit{Proof: }To prove this theorem, we refer to $[12]$.

\section{\protect\bigskip Fourier Transform and Plancherel Theorem on $AN.$}

\bigskip \textbf{4.1.} Let $G=SL(n,\mathbb{R})$\ be the real semi-simple Lie
group and let $G=KAN$ \ be the Iwasawa decomposition of $G$, where $K=SO(n,%
\mathbb{R}),$and%
\begin{equation}
N=\left( 
\begin{array}{ccccc}
1 & \ast & . & . & \ast \\ 
0 & 1 & \ast & . & \ast \\ 
. & . & . & . & \ast \\ 
. & . & . & . & . \\ 
0 & 0 & . & 0 & 1%
\end{array}%
\right) ,
\end{equation}

\begin{equation}
A=\left( 
\begin{array}{ccccc}
a_{1} & 0 & 0 & . & 0 \\ 
0 & a_{2} & 0 & . & 0 \\ 
. & . & . & . & . \\ 
. & . & . & . & . \\ 
0 & 0 & . & 0 & a_{n}%
\end{array}%
\right)
\end{equation}%
where $a_{1}.a_{2}....a_{n}=1$ and $a_{i}\in \mathbb{R}_{+}^{\ast }.$ The
product $AN$ is a closed subgroup of $G$ and is isomorphic (algebraically
and topologically) to the semi-direct product of $A$ and $N$ with $N$ normal
in $AN.$

Then the group $AN$ is nothing but the group $S=$ $N\rtimes $ $_{\rho }A.$
So the product of two elements $X$ and$Y$ by%
\begin{equation}
(x,\text{ }a)(m,\text{ }b)=(x.\rho (a)y,\text{ }a.b)
\end{equation}%
for\ any $X=(x,a_{1},a_{2},..,a_{n})\in S$ and $Y$ $%
=(m,b_{1},b_{2},..,b_{n})\in S.$ Let $dnda=dmda_{1}da_{2}..da_{n-1}$ be the
right haar measure on $S$ and let $L^{2}(S)$ be the Hilbert space of the
group $S.$ Let $L^{1}(S)$ be the Banach algebra that consists of all complex
valued functions on the group $S$, which are integrable with respect to the
Haar measure of $S$ and multiplication is defined by convolution on $S$ as

\begin{equation}
g\ast f=\int\limits_{S}f((m,b)^{-1}(n,a))g(m,b)dmdb
\end{equation}%
where $dmdb=dmdb_{1}db_{2}..db_{n-1}$ is the right Haar measure on $S=$ $%
N\rtimes $ $_{\rho }A.$

\bigskip In the following we prove the Plancherel theorem. Therefore let $%
T=N\times A$ be the Lie group of direct product of the two Lie groups $N$
and $A,$ and let $H=N\times A\times A$ the Lie group, with law%
\begin{equation}
(n,t,r)(m,s,q)=(n\rho (r)m,ts,rq)
\end{equation}%
for all $(n,t,r)\in H$ and $(m,s,q)\in H.$ In this case the group $S$ can be
identified with the closed subgroup $N\times \left\{ 0\right\} \times _{\rho
}A$ of $H$ and $T$ with the subgroup $N\times A\times \left\{ 0\right\} $of $%
H$

\ \textbf{Definition 4.1.}\textit{\ For every function }$f$ defined on $S$, 
\textit{one can define a function} $\widetilde{f}$ on $L$ \textit{as follows:%
} 
\begin{equation}
\widetilde{f}(n,a,b)=f(\rho (a)n,ab)
\end{equation}%
\textit{for all} $(n,a,b)\in H.$ \textit{So every function} $\psi (n,a)$ 
\textit{on} $S$\textit{\ extends uniquely as an invariant function} $%
\widetilde{\psi }(n,$ $b,$ $a)$ \textit{on} $L.$

\ \textbf{Remark 4.1. }\textit{The function} $\widetilde{f}$ \textit{is
invariant in the following sense:} 
\begin{equation}
\widetilde{f}(\rho (s)n,as^{-1},bs)=\widetilde{f}(n,a,b)
\end{equation}%
\textit{for any} $(n,a,b)\in H$ \textit{and} $s\in H.$

\textbf{Lemma 4.1.}\textit{\ For every function} $f\in L^{1}(S)$ \textit{and
for every }$g\mathcal{\in }$ $L^{1}(S)$, \textit{we have} 
\begin{equation}
g\ast \widetilde{f}(n,a,b)=g\ast _{c}\widetilde{f}(n,a,b)
\end{equation}%
\begin{equation}
\int\limits_{\mathbb{R}^{n-1}}\mathcal{F}_{A}^{n-1}\mathcal{F}^{d}\mathcal{(}%
g\ast \widetilde{f})(\lambda ,\mu ,\nu )d\nu =\mathcal{F}_{A}^{n-1}\mathcal{F%
}^{d}\widetilde{f}(\lambda ,\mu ,0)\mathcal{F}_{A}^{n-1}\mathcal{F}%
^{d}g(\lambda ,\mu )
\end{equation}%
\textit{for every} $(n,a,b)$ $\in H$, \textit{where} $\ast $\textit{\
signifies the convolution product on} $S$\ \textit{with respect the variables%
} $(n,b)$\ \textit{and} $\ast _{c}$\textit{signifies the commutative
convolution product on} $B$ \textit{with respect the variables }$(n,a),$ 
\textit{where }$\mathcal{F}_{A}^{n-1}$ is the Fourier transform on $A.$

\textbf{Plancherel} \textbf{Theorem 4.1. }\textit{For any function }$\Psi
\in L^{1}(S),$ \textit{we have \ }%
\begin{equation}
\int\limits_{\mathbb{R}^{d}}\int\limits_{\mathbb{R}^{n-1}}\left\vert 
\mathcal{F}^{d}\mathcal{F}_{A}^{n-1}\Psi (\lambda ,\mu )\right\vert
^{2}d\lambda =\int\limits_{AN}\left\vert \Psi (X,a)\right\vert ^{2}\text{ }dX
\end{equation}

\textit{Proof}: To prove this theorem, we refer to $[12]$.

\section{\protect\bigskip Fourier Transform and Plancherel Theorem on $SL(n,%
\mathbb{R}).$}

\bigskip \textbf{5.1.} In the following we use the Iwasawa decomposition of $%
G=SL(n,\mathbb{R}),$ to define the Fourier transform and to get Plancherel
theorem on $G=SL(n,\mathbb{R}).$ We denote by $L^{1}(G)$ the Banach algebra
that consists of all complex valued functions on the group $G$, which are
integrable with respect to the Haar measure of $G$ and multiplication is
defined by convolution on $G$ , 
\begin{equation}
\phi \ast f(g)=\int\limits_{G}f(h^{-1}g)\phi (g)dg
\end{equation}

Let $G=SL(n,\mathbb{R})=KNA$ be the Iawsawa decomposition of $G.$The Haar
measure $dg$ on $G$ can be calculated from the Haar measures $dn;$ $da$ and $%
dk$ on $N;A$ and $K;$ respectively, by the formula%
\begin{equation}
\int\limits_{G}f(g)dg=\int\limits_{A}\int\limits_{N}\int%
\limits_{K}f(ank)dadndk
\end{equation}

Keeping in mind that $a^{-2\rho }$ is the modulus of the automorphism $%
n\rightarrow $ $ana^{-1}$ of $N$ we get also the following representation of 
$dg$ 
\begin{equation}
\int\limits_{G}f(g)dg=\int\limits_{A}\int\limits_{N}\int%
\limits_{K}f(ank)dadndk=\int\limits_{N}\int\limits_{A}\int%
\limits_{K}f(nak)a^{-2^{\rho }}dndadk
\end{equation}%
where $\rho =\dim N$ $=\frac{n(n-1)}{2}=1+2+3+....+n-2+n-1.$ Furthermore,
using the relation $\int\limits_{G}f(g)dg=\int\limits_{G}f(g^{-1})dg,$ we
receive 
\begin{equation}
\int\limits_{K}\int\limits_{A}\int\limits_{N}f(nak)a^{-2\rho
}dndadk=\int\limits_{K}\int\limits_{A}\int\limits_{N}f(kan)a^{2\rho }dndadk
\end{equation}

Let $\underline{k}$ be the Lie algebra of $K$ and $(X_{1},X_{2},.....,X_{m})$
a basis of $\underline{k}$ , such that the both operators%
\begin{equation}
\Delta =\dsum\limits_{i=1}^{m}X_{i}^{2}
\end{equation}%
\begin{equation}
D_{q}=\dsum\limits_{0\leq l\leq q}\left(
-\dsum\limits_{i=1}^{m}X_{i}^{2}\right) ^{l}
\end{equation}%
are left and right invariant (bi-invariant) on $K,$ this basis exist see $%
[2, $ $p.564)$. For $l\in 
\mathbb{N}
$, let $D^{l}=(1-\Delta )^{l}$, then the family of semi-norms $\{\sigma _{l}$%
, $l\in 
\mathbb{N}
\}$ such that%
\begin{equation}
\sigma _{l}(f)=(\int_{K}\left\vert D^{l}f(y)\right\vert ^{2}dy)^{\frac{1}{2}%
},\text{ \ \ \ \ \ \ \ \ }f\in C^{\infty }(K)
\end{equation}%
define on $C^{\infty }(K)$ the same topology of the Frechet topology defined
by the semi-normas $\left\Vert X^{\alpha }f\right\Vert _{2}$ defined as%
\begin{equation}
\left\Vert X^{\alpha }f\right\Vert _{2}=(\int_{K}\left\vert X^{\alpha
}f(y)\right\vert ^{2}dy)^{\frac{1}{2}},\text{ \ \ \ \ \ \ \ \ }f\in
C^{\infty }(K)
\end{equation}%
where $\alpha =(\alpha _{1},$.....,$\alpha _{m})\in 
\mathbb{N}
^{m},$ for all the above formula see $[2,P.176-177]$ and $[2,p.565]$

Let $\widehat{K}$ be the set of all irreducible unitary representations of $%
K.$ If $\gamma \in \widehat{K}$, we denote by $E_{\gamma }$ the space of
representation $\gamma $ and $d_{\gamma }$\ its dimension \qquad

\textbf{Definition 5.1.} \textit{The Fourier transform of a function }$f\in
C^{\infty }(K)$\textit{\ is defined as} 
\begin{equation}
Tf(\gamma )=\dint\limits_{K}f(x)\gamma (x^{-1})dx
\end{equation}%
\textit{where }$T$\textit{\ is the Fourier transform on} $K$

\textbf{Theorem (A. Cerezo) 5.1.} \textit{Let} $f\in C^{\infty }(K),$ 
\textit{then we have the inversion of the Fourier transform} 
\begin{equation}
f(x)=\dsum\limits_{\gamma \in \widehat{K}}d\gamma tr[Tf(\gamma )\gamma (x)]
\end{equation}%
\begin{equation}
f(x^{-1})=\dsum\limits_{\gamma \in \widehat{K}}d\gamma tr[Tf(\gamma )\gamma
(x^{-1})]=\dsum\limits_{\gamma \in \widehat{K}}d\gamma tr[Tf(\gamma )\gamma
^{\ast }(x)]
\end{equation}

\begin{equation}
f(I_{K})=\dsum\limits_{\gamma \in \widehat{K}}d\gamma tr[Tf(\gamma )]
\end{equation}%
\textit{and the Plancherel formula} 
\begin{equation}
\left\Vert f(x)\right\Vert _{2}^{2}=\dint \left\vert f(x)\right\vert
^{2}dx=\dsum\limits_{\gamma \in \widehat{K}}d_{\gamma }\left\Vert Tf(\gamma
)\right\Vert _{H.S}^{2}
\end{equation}%
\textit{for any }$f\in L^{1}(K),$ \textit{where }$I_{K}$ \textit{is the
identity element of \ }$K,$ see $[2,P.562-563],$\textit{where }$\left\Vert
Tf(\gamma )\right\Vert _{H.S}^{2}$ \textit{is} the \textit{norm of
Hilbert-Schmidt of the operator }$Tf(\gamma ).$

\bigskip \textbf{Definition 5.2}. \textit{For any function} $f\in \mathcal{D}%
(G),$ \textit{we can define a function} $\Upsilon (f)$\textit{on }$G\times K$
\textit{by} 
\begin{equation}
\Upsilon (f)(g,k_{1})=\Upsilon (f)(kna,k_{1})=f(gk_{1})=f(knak_{1})
\end{equation}%
\textit{for }$g=kna\in G,$ \textit{and} $k_{1}\in K$ . \textit{The
restriction of} $\ \Upsilon (f)\ast \psi (g,k_{1})$ \textit{on} $K(G)$ 
\textit{is }$\Upsilon (f)\ast \psi (g,k_{1})\downarrow _{K(G)}=f(nak_{1})\in 
\mathcal{D}(G),$ \textit{and }$\Upsilon (f)(g,k_{1})\downarrow
_{K}=f(g,I_{K})=f(kna)$ $\in \mathcal{D}(G)$

\textbf{Remark 5.1. }\textit{The function }$\Upsilon (f)$ \textit{is
invariant in the following sense }%
\begin{equation}
\Upsilon (f)(gh,h^{-1}k_{1})=f(gk_{1})=f(knak_{1})
\end{equation}

\textbf{Definition 5.3}.\textit{\ Let }$f$ \textit{and }$\psi $ \textit{be
two functions belong to} $\mathcal{D}(G),$ \textit{then we can define the
convolution of } $\Upsilon (f)\ $\textit{and} $\psi $\ \textit{on} $G$ $%
\times K$ \textit{as}

\begin{eqnarray*}
\Upsilon (f)\ast \psi (g,k_{1}) &=&\int\limits_{G}\Upsilon
(f)(gg_{2}^{-1},k_{1})\psi (g_{2})dg_{2} \\
&=&\int\limits_{K}\int\limits_{N}\int\limits_{A}\Upsilon
(f)(knaa_{2}^{-1}n_{2}^{-1}k^{-1}k_{1})\psi
(k_{2}n_{2}a_{2})dk_{2}dn_{2}da_{2}
\end{eqnarray*}

So we get 
\begin{eqnarray*}
\Upsilon (f)\ast \psi (g,k_{1}) &\downarrow &_{K(G)}=\Upsilon (f)\ast \psi
(I_{K}na,k_{1}) \\
&=&\int\limits_{K}\int\limits_{N}\int%
\limits_{A}f(naa_{2}^{-1}n_{2}^{-1}k^{-1}k_{1})\psi
(k_{2}n_{2}a_{2})dk_{2}dn_{2}da_{2}
\end{eqnarray*}%
where $T$ is the Fourier transform on $K,$ and\textit{\ }$I_{K}$ \textit{i}s
the identity element of\textit{\ }$K.$ Denote by\textit{\ }$\mathcal{F}$ is
the Fourier transform on $AN$

\textbf{Definition} \textbf{5.4.} \textit{If }$f\in \mathcal{D}(G)$ \textit{%
and let} $\Upsilon (f)$ \textit{be the associated function to} $f$ , \textit{%
we define the Fourier transform of \ }$\Upsilon (f)(g,k_{1})$ \textit{by }%
\begin{eqnarray}
&&T\mathcal{F}\Upsilon (f))(I_{K},\xi ,\lambda ,\gamma )=T\mathcal{F}%
\Upsilon (f))(I_{K},\xi ,\lambda ,\gamma )  \notag \\
&=&\int_{K}\int_{N}\int_{A}\dsum\limits_{\delta \in \widehat{K}}d_{\delta
}tr[\int_{K}\Upsilon (f)(kna,k_{1})\delta (k^{-1})dk]a^{-i\lambda }e^{-\text{
}i\langle \text{ }\xi ,\text{ }n\rangle }\text{ }\gamma
(k_{1}^{-1})dadndk_{1}  \notag \\
&=&\int_{N}\int_{A}\int_{K}\Upsilon (f)(I_{K}na,k_{1})a^{-i\lambda }e^{-%
\text{ }i\langle \text{ }\xi ,\text{ }n\rangle }\text{ }\gamma
(k_{1}^{-1})dadndk_{1}  \notag \\
&=&\int_{N}\int_{A}\int_{K}f(nak_{1})a^{-i\lambda }e^{-\text{ }i\langle 
\text{ }\xi ,\text{ }n\rangle }\text{ }\gamma (k_{1}^{-1})dadndk_{1}
\end{eqnarray}

\textbf{Theorem 5.2.} (\textbf{\textit{Plancherel's Formula for the Group }}$%
G$ )\textbf{\textit{. }}\textit{For any function\ }$f\in $\textit{\ }$%
L^{1}(G)\cap $\textit{\ }$L^{2}(G),$\textit{we get }%
\begin{eqnarray}
\int \left\vert f(g)\right\vert ^{2}dg &=&\int_{K}\int_{N}\int_{A}\left\vert
f(kna)\right\vert ^{2}dadndk  \notag \\
&=&\sum_{\gamma \in \widehat{K}}d_{\gamma }\int\limits_{\mathbb{R}%
^{d}}\int\limits_{\mathbb{R}^{n-1}}\left\Vert T\mathcal{F}f(\lambda ,\xi
,\gamma )\right\Vert _{H.S}^{2}d\lambda d\xi  \notag \\
&=&\int\limits_{\mathbb{R}^{d}}\int\limits_{\mathbb{R}^{n-1}}\sum_{\gamma
\in \widehat{K}}d_{\gamma }\left\Vert T\mathcal{F}f(\lambda ,\xi ,\gamma
)\right\Vert _{H.S}^{2}d\lambda d\xi
\end{eqnarray}%
\textit{and so if }$f\in L^{1}(G),$ \textit{then the Fourier inversion is\ }%
\begin{eqnarray}
f(ank_{1}) &=&\int\limits_{\mathbb{R}^{d}}\int\limits_{\mathbb{R}%
^{n-1}}\sum_{\gamma \in \widehat{K}}d_{\gamma }tr[T\mathcal{F}f((\lambda
,\xi ,\gamma )\gamma (k_{1})]a^{i\lambda }e^{\text{ }i\langle \text{ }\xi ,%
\text{ }n\rangle }d\lambda d\xi  \notag \\
f(I_{A}I_{N}I_{K}) &=&\int\limits_{\mathbb{R}^{d}}\int\limits_{\mathbb{R}%
^{n-1}}\sum_{\gamma \in \widehat{K}}d_{\gamma }tr[T\mathcal{F}f(\lambda ,\xi
,\gamma )]d\lambda d\xi
\end{eqnarray}%
\textit{where} $I_{A},I_{N}$ and $I_{K}$ \textit{are the identity elements of%
} $A$, $N$ \textit{and }$K$ \textit{respectively, where }$\mathcal{F}$ 
\textit{is the Fourier transform on }$AN$ \textit{and }$T$ \textit{is the
Fourier transform on} $K,$ \textit{and }$I_{K}$ \textit{is the identity
element of }$K$

\bigskip \textit{Proof: }First let $\overset{\vee }{f}$ be the function
defined by 
\begin{equation}
\ \overset{\vee }{f}(kna)=\overline{f((kna)^{-1})}=\overline{%
f(a^{-1}n^{-1}k^{-1})}
\end{equation}

Then we have%
\begin{eqnarray}
&&\int \left\vert f(g)\right\vert ^{2}dg  \notag \\
&=&\Upsilon (f)\ast \overset{\vee }{f}(I_{K}I_{N}I_{A},I_{K_{1}})  \notag \\
&=&\int\limits_{G}\Upsilon (f)(I_{K}I_{N}I_{A}(g_{2}^{-1}),I_{K_{1}})\overset%
{\vee }{f}(g_{2})dg_{2}  \notag \\
&=&\int\limits_{A}\int\limits_{N}\int\limits_{K}\Upsilon
(f)(a_{2}^{-1}n_{2}^{-1}k_{2}^{-1},I_{K})\overset{\vee }{f}%
(k_{2}n_{2}a_{2})da_{2}dn_{2}dk_{2}  \notag \\
&=&\int\limits_{A}\int\limits_{N}\int%
\limits_{K}f(a_{2}^{-1}n_{2}^{-1}k_{2}^{-1})\overline{%
f((k_{2}n_{2}a_{2})^{-1})}da_{2}dn_{2}dk_{2}  \notag \\
&=&\int\limits_{A}\int\limits_{N}\int\limits_{K}\left\vert
f(a_{2}n_{2}k_{2})\right\vert ^{2}da_{2}dn_{2}dk_{2}
\end{eqnarray}

\bigskip In other hand%
\begin{eqnarray*}
&&\Upsilon (f)\ast \overset{\vee }{f}(I_{K}I_{N}I_{A},I_{K_{1}}) \\
&=&\int\limits_{\mathbb{R}^{d}}\int\limits_{\mathbb{R}^{n-1}}\text{ }%
\mathcal{F}(\Upsilon (f)\ast \overset{\vee }{f})(I_{K},\lambda ,\xi
,I_{K_{1}})d\lambda d\xi \\
&=&\int\limits_{\mathbb{R}^{d}}\int\limits_{\mathbb{R}^{n-1}}\sum_{\gamma
\in \widehat{K}}d_{\gamma }\sum_{\delta \in \widehat{K}}d_{\delta }tr[T%
\mathcal{F}(\Upsilon (f)\ast \overset{\vee }{f})(\delta ,\lambda ,\xi
,\gamma )]d\lambda d\xi \\
&=&\int\limits_{\mathbb{R}^{d}}\int\limits_{\mathbb{R}^{n-1}}\sum_{\gamma
\in \widehat{K}}d_{\gamma }tr[\dint\limits_{K}\mathcal{F}(\Upsilon (f)\ast 
\overset{\vee }{f})(I_{K},\lambda ,\xi ,k_{1})\gamma
(k_{1}^{-1})dk_{1}]d\lambda d\xi \\
&=&\int\limits_{\mathbb{R}^{d}}\int\limits_{\mathbb{R}^{n-1}}\dint%
\limits_{N}\dint\limits_{A}\sum_{\gamma \in \widehat{K}}d_{\gamma
}tr[\dint\limits_{K}\mathcal{F}(\Upsilon (f)\ast \overset{\vee }{f}%
)(I_{k}na,k_{1})\gamma (k_{1}^{-1})dk_{1}] \\
&&a^{-i\lambda }e^{-\text{ }i\langle \text{ }\xi ,\text{ }n\rangle
}dndad\lambda d\xi \\
&=&\int\limits_{\mathbb{R}^{d}}\int\limits_{\mathbb{R}^{n-1}}\dint%
\limits_{N}\dint\limits_{A}\dint\limits_{A}\int\limits_{N}\sum_{\gamma \in 
\widehat{K}}d_{\gamma }tr[\dint\limits_{K}\dint\limits_{K}\Upsilon
(f)(I_{k}naa_{2}^{-1}n_{2}^{-1}k_{2}^{-1},k_{1})\overset{\vee }{f}%
(k_{2}n_{2}a_{2})\gamma (k_{1}^{-1})dk_{1}dk_{2}] \\
&&a^{-i\lambda }e^{-\text{ }i\langle \text{ }\xi ,\text{ }n\rangle
}dndada_{2}dn_{2}d\lambda d\xi \\
&=&\int\limits_{\mathbb{R}^{d}}\int\limits_{\mathbb{R}^{n-1}}\dint%
\limits_{N}\dint\limits_{A}\dint\limits_{A}\int\limits_{N}\sum_{\gamma \in 
\widehat{K}}d_{\gamma }tr[\dint\limits_{K}\dint\limits_{K}\Upsilon
(f)(aa_{2}^{-1}nn_{2}^{-1}k_{2}^{-1},k_{1})\overset{\vee }{f}%
(k_{2}n_{2}a_{2})\gamma (k_{1}^{-1})dk_{1}dk_{2}] \\
&&a^{-i\lambda }e^{-\text{ }i\langle \text{ }\xi ,\text{ }n\rangle
}dndada_{2}dn_{2}d\lambda d\xi \\
&=&\int\limits_{\mathbb{R}^{d}}\int\limits_{\mathbb{R}^{n-1}}\dint%
\limits_{N}\dint\limits_{A}\dint\limits_{A}\int\limits_{N}\sum_{\gamma \in 
\widehat{K}}d_{\gamma
}tr[\dint\limits_{K}\dint\limits_{K}f(aa_{2}^{-1}nn_{2}^{-1}k_{2}^{-1}k_{1})%
\overset{\vee }{f}(k_{2}n_{2}a_{2})\gamma (k_{1}^{-1})dk_{1}dk_{2}] \\
&&a^{-i\lambda }e^{-\text{ }i\langle \text{ }\xi ,\text{ }n\rangle
}dndada_{2}dn_{2}d\lambda d\xi \\
&=&\int\limits_{\mathbb{R}^{d}}\int\limits_{\mathbb{R}^{n-1}}\dint%
\limits_{N}\dint\limits_{A}\dint\limits_{A}\int\limits_{N}\sum_{\gamma \in 
\widehat{K}}d_{\gamma
}tr[\dint\limits_{K}\dint\limits_{K}f(aa_{2}^{-1}nn_{2}^{-1}k_{1})\overset{%
\vee }{f}(k_{2}n_{2}a_{2})\gamma (k_{2}^{-1})\gamma (k_{1}^{-1})dk_{1}dk_{2}]
\\
&&a^{-i\lambda }e^{-\text{ }i\langle \text{ }\xi ,\text{ }n\rangle
}dndada_{2}dn_{2}d\lambda d\xi
\end{eqnarray*}

\bigskip We have used the fact that%
\begin{equation}
\int\limits_{A}\int\limits_{N}\int\limits_{K}f(kna)dadndk=\int\limits_{N}%
\int\limits_{A}\int\limits_{K}f(kan)a^{2\rho }dndadk
\end{equation}%
and 
\begin{eqnarray}
&&\int\limits_{\mathbb{R}^{d}}\int\limits_{A}\int\limits_{N}\int%
\limits_{K}f(kna)e^{-\text{ }i\langle \text{ }(\xi ,\text{ }n\rangle
}dadndkd\xi  \notag \\
&=&\int\limits_{\mathbb{R}^{d}}\int\limits_{A}\int\limits_{N}\int%
\limits_{K}f(kan)e^{-\text{ }i\langle \text{ }\xi ,\text{ }ana^{-1}\text{ }%
\rangle }a^{2\rho }dadndkd\xi  \notag \\
&=&\int\limits_{\mathbb{R}^{d}}\int\limits_{A}\int\limits_{N}\int%
\limits_{K}f(kan)e^{-\text{ }i\langle \text{ }a\xi a^{-1},\text{ }n\rangle
}a^{2\rho }dadndkd\xi  \notag \\
&=&\int\limits_{\mathbb{R}^{d}}\int\limits_{A}\int\limits_{N}\int%
\limits_{K}f(kan)e^{-\text{ }i\langle \text{ }\xi ,\text{ }n\rangle
}dadndkd\xi
\end{eqnarray}%
where $\langle \xi ,n\rangle =\dsum\limits_{i=1}^{d}\xi _{i}n_{i},$ $d\xi
=d\xi _{1}d\xi _{2}$, and $ana^{-1}=a(\xi )a^{-1}=a^{-2\rho }\xi $, then we
get%
\begin{eqnarray*}
&&\Upsilon \mathcal{(}f)\ast \overset{\vee }{f}(I_{K}I_{N}I_{A},I_{K_{1}}) \\
&=&\int\limits_{\mathbb{R}^{d}}\int\limits_{\mathbb{R}^{n-1}}\dint%
\limits_{A}\dint\limits_{N}\dint\limits_{A}\int\limits_{N}\sum_{\gamma \in 
\widehat{K}}d_{\gamma }tr[\dint\limits_{K}\dint\limits_{K}f(ank_{1})\overset{%
\vee }{f}(k_{2}n_{2}a_{2})\gamma (k_{2}^{-1})\gamma (k_{1}^{-1})dk_{1}dk_{2}]
\\
&&a^{-i\lambda }e^{-\text{ }i\langle \text{ }\xi ,\text{ }n\rangle
}a_{2}^{-i\lambda }e^{-\text{ }i\langle \text{ }\xi ,\text{ }n_{2}\rangle
}dada_{2}dndn_{2}d\lambda d\xi \\
&=&\int\limits_{\mathbb{R}^{d}}\int\limits_{\mathbb{R}^{n-1}}\dint%
\limits_{A}\dint\limits_{N}\dint\limits_{A}\int\limits_{N}\sum_{\gamma \in 
\widehat{K}}d_{\gamma }tr[\dint\limits_{K}\dint\limits_{K}f(ank_{1})%
\overline{f((k_{2}n_{2}a_{2})^{-1})}\gamma (k_{2}^{-1})\gamma
(k_{1}^{-1})dk_{1}dk_{2}] \\
&&a^{-i\lambda }e^{-\text{ }i\langle \text{ }\xi ,\text{ }n\rangle
}a_{2}^{-i\lambda }e^{-\text{ }i\langle \text{ }\xi ,\text{ }n_{2}\rangle
}dndadn_{2}da_{2}d\lambda d\xi \\
&=&\int\limits_{\mathbb{R}^{d}}\int\limits_{\mathbb{R}^{n-1}}\dint%
\limits_{A}\dint\limits_{N}\dint\limits_{A}\int\limits_{N}\sum_{\gamma \in 
\widehat{K}}d_{\gamma }tr[\dint\limits_{K}\dint\limits_{K}f(ank_{1})%
\overline{f((a_{2}{}^{-1}n_{2}^{-1}k_{2}^{-1})}\gamma (k_{2}^{-1})\gamma
(k_{1}^{-1})dk_{1}dk_{2}] \\
&&a^{-i\lambda }e^{-\text{ }i\langle \text{ }\xi ,\text{ }n\rangle
}a_{2}^{-i\lambda }e^{-\text{ }i\langle \text{ }\xi ,\text{ }n_{2}\rangle
}dndadn_{2}da_{2}d\lambda d\xi \\
&=&\int\limits_{\mathbb{R}^{d}}\int\limits_{\mathbb{R}^{n-1}}\dint%
\limits_{A}\dint\limits_{N}\dint\limits_{A}\int\limits_{N}\sum_{\gamma \in 
\widehat{K}}d_{\gamma }tr[\dint\limits_{K}\dint\limits_{K}f(ank_{1})%
\overline{f((a_{2}{}n_{2}k_{2})}\gamma ^{\ast }(k_{2}^{-1})\gamma
(k_{1}^{-1})dk_{1}dk_{2}] \\
&&a^{-i\lambda }e^{-\text{ }i\langle \text{ }\xi ,\text{ }n\rangle }%
\overline{a_{2}^{-i\lambda }}\text{ }\overline{e^{\text{ -}i\langle \text{ }%
\xi ,\text{ }n_{2}\rangle }}dndadn_{2}da_{2}d\lambda d\xi \\
&=&\int\limits_{\mathbb{R}^{d}}\int\limits_{\mathbb{R}^{n-1}}\sum_{\gamma
\in \widehat{K}}d_{\gamma }T\mathcal{F}f(\lambda ,\xi ,\gamma )\overline{T%
\mathcal{F}f(\lambda ,\xi ,\gamma ^{\ast })}d\lambda d\xi \\
&=&\int\limits_{\mathbb{R}^{d}}\int\limits_{\mathbb{R}^{n-1}}\sum_{\gamma
\in \widehat{K}}d_{\gamma }\left\Vert T\mathcal{F}(f)(\lambda ,\xi ,\gamma
)\right\Vert _{H.S}^{2}d\lambda d\xi
\end{eqnarray*}

\section{\protect\bigskip Fourier Transform and Plancherel Theorem on $GL(n,%
\mathbb{R)}.$}

\textbf{6.1. New Group. }Let $GL(n,\mathbb{R)}$ be the general linear group
consisting of all matrices of the form 
\begin{equation}
GL(n,\mathbb{R)}=\{X=\left( 
\begin{array}{c}
a_{ij}%
\end{array}%
\right) ,a_{ij}\in \mathbb{R},\text{ }\prec 1\prec n,,\text{ }\prec j\text{ }%
\prec n,\text{ }and\text{ }\det A\neq 0\}
\end{equation}

As a manifold, $GL(n,\mathbb{R})$ is not connected but rather has two
connected components: the matrices with positive determinant and the ones
with negative determinant which is denoted by $GL_{-}(n,\mathbb{R})$. The
identity component, denoted by $GL_{+}(n,\mathbb{R})$, consists of the real $%
n\times n$ matrices with positive determinant. This is also a Lie group of
dimension $n^{2}$; it has the same Lie algebra as $GL(n,\mathbb{R})$.

The group $GL(n,\mathbb{R})$ is also noncompact. The maximal compact
subgroup of $GL(n,\mathbb{R})$ is the orthogonal group $O(n)$, while the
maximal compact subgroup of $GL_{+}(n,\mathbb{R})$ is the special orthogonal
group $SO(n)$. As for $SO(n)$, the group $GL_{+}(n,\mathbb{R})$ is not
simply connected

\textbf{Theorem 6.1. }$GL_{-}(n,\mathbb{R)}$ \textit{is group isomorphic
onto }$GL_{+}(n,\mathbb{R})$

\bigskip \textit{Proof:} $GL_{-}(n,\mathbb{R)}$ is the subset of $GL(n,%
\mathbb{R)},$ which is defined as 
\begin{equation}
GL_{-}(n,\mathbb{R)}=\{A\in GL(n,\mathbb{R)},\text{ }\det A\text{ }\langle 0%
\text{ }\}
\end{equation}%
We supply $GL_{-}(n,\mathbb{R)}$ by the law noted\textit{\ }$\bullet $ which
is defined by\textit{\ }%
\begin{equation}
A\bullet B=I^{-}A.B
\end{equation}%
for any $A\in GL_{-}(n,\mathbb{R)}$ and $B\in GL_{-}(n,\mathbb{R)}$, where $%
. $ signifies the usual multiplication of two matrix in $GL_{-}(n,\mathbb{R)}
$ and $I_{-}$ is the the following matrix defined as%
\begin{equation}
A=\left( \left[ a_{ij}\right] \right)
\end{equation}%
and%
\begin{equation}
I^{-}=\left( 
\begin{array}{ccccc}
-1 & 0 & 0 & . & 0 \\ 
0 & 1 & 0 & . & 0 \\ 
. & . & . & . & . \\ 
. & . & . & . & . \\ 
0 & 0 & . & 0 & 1%
\end{array}%
\right)
\end{equation}

Then we have $A\bullet B$ $\in GL_{-}(n,\mathbb{R)}$, for any $A\in GL_{-}(n,%
\mathbb{R)}$ and $B\in GL_{-}(n,\mathbb{R)}$.

The identity element is $I^{-}$ because if

\begin{equation}
A\bullet B=I^{-}.A.B=B
\end{equation}%
then we have $A=I^{-}$ and so 
\begin{equation}
I^{-}\bullet B=B\bullet I^{-}
\end{equation}

The law is associative, because%
\begin{eqnarray*}
A\bullet (B\bullet C) &=&I^{-}.A.(B\bullet C)=I^{-}.A.(I^{-}.B.C) \\
&=&A.B.C=(I^{-}.A.B).(I^{-}.C)=(A\bullet B)(I^{-}.C) \\
&=&I^{-}.(A\bullet B).C=(A\bullet B)\bullet C
\end{eqnarray*}%
and it is easy th show the inverse of any element $A\in GL_{-}(n,\mathbb{R)}$
is

\begin{equation}
I^{-}.A
\end{equation}

\bigskip Now consider the mapping $\varphi :A\longrightarrow B$ defined by%
\begin{equation}
\varphi (A)=I^{-}.A
\end{equation}%
for any $A\in GL_{-}(n,\mathbb{R)}$. Then we get%
\begin{eqnarray}
\varphi (A\bullet B) &=&I^{-}(A\bullet B)=I^{-}.I^{-}.A.B  \notag \\
&=&I^{-}.A.I^{-}.B=\varphi (A).\varphi (B)
\end{eqnarray}

\bigskip It is obvious that $\varphi $ is one-to-one and onto, so $\varphi $
is a group isomorphism from $GL_{-}(n,\mathbb{R)}$\ onto $GL_{+}(n,\mathbb{R)%
}.$ As well known the group $GL_{+}(n,\mathbb{R)}$ is isomorphic onto the
direct product of the two groups $SL(n,\mathbb{R)}$ and $\mathbb{R}%
_{+}^{\ast },$ $i.e$ $GL_{+}(n,\mathbb{R)=}SL(n,\mathbb{R)\times }$ $\mathbb{%
R}_{+}^{\ast }$ and $GL(n,\mathbb{R)=}GL_{-}(n,\mathbb{R)}\ \cup GL_{+}(n,%
\mathbb{R)=(}SL(n,\mathbb{R)\times }$ $\mathbb{R}_{+}^{\ast })\cup \mathbb{(}%
SL(n,\mathbb{R)\times }$ $\mathbb{R}_{+}^{\ast }).$ Our aim result is

\textbf{Plancherel theorem 6.2. }\textit{Let} $\mathcal{F}_{+}^{\ast }$ 
\textit{be the Fourier transform on }$GL_{+}(n,\mathbb{R}$ \textit{the we get%
} 
\begin{eqnarray}
\int_{GL_{+}(n,\mathbb{R)}}\left\vert f(g,t)\right\vert ^{2}dg\frac{dt}{t}
&=&\int\limits_{\mathbb{R}_{+}^{\ast }}\int_{K}\int_{N}\int_{A}\left\vert
f(kna,t)\right\vert ^{2}dadndk\frac{dt}{t}  \notag \\
&=&\sum_{\gamma \in \widehat{K}}d_{\gamma }\int\limits_{\mathbb{R}%
}\int\limits_{\mathbb{R}^{d}}\int\limits_{\mathbb{R}^{n-1}}\left\Vert 
\mathcal{F}_{+}^{\ast }T\mathcal{F}f(\lambda ,\xi ,\gamma ,\eta )\right\Vert
_{H.S}^{2}d\lambda d\xi d\eta  \notag \\
&=&\int\limits_{\mathbb{R}}\int\limits_{\mathbb{R}^{d}}\int\limits_{\mathbb{R%
}^{n-1}}\sum_{\gamma \in \widehat{K}}d_{\gamma }\left\Vert \mathcal{F}%
_{+}^{\ast }T\mathcal{F}f(\lambda ,\xi ,\gamma ,\eta )\right\Vert
_{H.S}^{2}d\lambda d\xi d\eta
\end{eqnarray}%
for any $f\in L^{1}(GL_{+}(n,\mathbb{R)}\cap L^{2}(GL_{+}(n,\mathbb{R)}$.
The proof of theorem results immediately from theorem 5.2

\textbf{Corollary 6.1. }\textit{Let }$f$ be a function belongs $L^{1}(GL(n,%
\mathbb{R)}\cap L^{2}(GL(n,\mathbb{R)}$ 
\begin{eqnarray}
&&\int_{GL(n,\mathbb{R)}}\left\vert f(g,t)\right\vert ^{2}dg\frac{dt}{t} 
\notag \\
&=&\int_{GL_{-}(n,\mathbb{R)}\cup GL_{+}(n,\mathbb{R)}}\left\vert
f(g,t)\right\vert ^{2}dg\frac{dt}{t}=2\int_{GL_{+}(n,\mathbb{R)}}\left\vert
f(g,t)\right\vert ^{2}dg\frac{dt}{t}  \notag \\
&=&2\int\limits_{\mathbb{R}_{+}^{\ast }}\int_{K}\int_{N}\int_{A}\left\vert
f(kna,t)\right\vert ^{2}dadndk\frac{dt}{t}  \notag \\
&=&2\sum_{\gamma \in \widehat{K}}d_{\gamma }\int\limits_{\mathbb{R}%
}\int\limits_{\mathbb{R}^{d}}\int\limits_{\mathbb{R}^{n-1}}\left\Vert 
\mathcal{F}_{+}^{\ast }T\mathcal{F}f(\lambda ,\xi ,\gamma ,\eta )\right\Vert
_{H.S}^{2}d\lambda d\xi d\eta  \notag \\
&=&2\int\limits_{\mathbb{R}}\int\limits_{\mathbb{R}^{d}}\int\limits_{\mathbb{%
R}^{n-1}}\sum_{\gamma \in \widehat{K}}d_{\gamma }\left\Vert \mathcal{F}%
_{+}^{\ast }T\mathcal{F}f(\lambda ,\xi ,\gamma ,\eta )\right\Vert
_{H.S}^{2}d\lambda d\xi d\eta
\end{eqnarray}

\textbf{Remark 6.1.} this corollary explains the Fourier transform and
Plancherel formula on the non connected Lie group $GL(n,\mathbb{R)}$.

\end{document}